\documentclass[11pt]{amsart}  
\DeclareMathOperator{\ad}{ad}
\DeclareMathOperator{\ch}{char}

\newtheorem{theorem}{Theorem}[section]

\newtheorem{lemma}[theorem]{Lemma}

\newtheorem{proposition}[theorem]{Proposition}

\theoremstyle{definition}

\newtheorem{example}[theorem]{Example}

\theoremstyle{remark}

\makeatletter
\def\romenumi{%
  \def\theenumi{\roman{enumi}}%
  \def\p@enumi{\theenumi}%
  \def\labelenumi{(\@roman\c@enumi)}}
\makeatother

\begin{document}  
\title{Associative algebras satisfying a semigroup identity} 
 
\author{D.~M.~Riley}  
\address{Department of Mathematics,  The University of Alabama, Tuscaloosa, AL  
35487-0350, USA}  
\email{driley@gp.as.ua.edu} 
 
\author{Mark~C.~Wilson}
\address{Department of Mathematics, University of Auckland, Private Bag 92019  
Auckland, New Zealand}  
\email{wilson@math.auckland.ac.nz}

\date{\today}  

\begin{abstract}  Denote by $(R,\cdot)$ the multiplicative semigroup of 
an associative algebra $R$ over an infinite field, and let
$(R,\circ)$ represent $R$ when viewed as a semigroup via the 
circle operation $x\circ y=x+y+xy$.  
In this paper we characterize the existence of an identity in
these semigroups in terms of the Lie structure of $R$.
Namely, we prove that the following conditions on $R$ are equivalent:
the semigroup $(R,\circ)$ satisfies an identity; 
the semigroup $(R,\cdot)$ satisfies a reduced identity; and,
the associated Lie algebra of $R$ satisfies the Engel condition.
When $R$ is finitely generated these conditions 
are each equivalent to $R$ being upper Lie nilpotent.
\end{abstract}

\maketitle
\section{Introduction and statement of results}

A well-known result due to Levitzki (\cite{levitzki:kurosch})  
states that every finitely generated bounded nil ring is nilpotent. 
Not long ago, Zel'manov proved the Lie-theoretic analogue: every 
finitely generated Lie ring satisfying the Engel condition is nilpotent 
(\cite{zelmanov:burnside}).
The corresponding problem in the category of groups is the 
famous Burnside problem.  The construction by Adian and Novikov 
 of infinite finitely generated groups of finite exponent  
provided a negative solution to this problem (see \cite{adian}).

The Burnside problem has some natural generalizations.
For example, the problem of whether or not every Engel group is
locally nilpotent remains open (\cite{shalev:p-survey}).
Because every nilpotent group is known to satisfy a semigroup identity
(\cite{malcev, neumann-taylor}), a weaker version of this problem has also been posed: 
does every Engel group satisfy a semigroup identity (\cite[Problem 2.82]{kourovka})?
Even the following question (\cite{rhemtulla:private}) remains
open: can an Engel group contain a free  (noncommutative)
subsemigroup? 

Recently, the present authors settled the ring-theoretic 
analogues of these problems.   

Recall that $R$ satisfies the {\em Engel identity} of degree $n$ if and 
only if
$$e_n:=[x,\underbrace{y,y,\ldots,y}_n]$$ 
is identically zero in $R$; whereas, $R$ is said to be {\em upper Lie nilpotent}  if the
descending central series of associative
ideals $\{\gamma^i(R)\}$ in $R$ defined by $\gamma^1(R)=R,
\gamma^{i+1}(R)=\langle [\gamma^i(R),R]\rangle$ reaches zero in finitely many steps. 
In addition to the usual multiplicative semigroup, $(R,\cdot)$, 
$R$ forms a semigroup, denoted by $(R,\circ)$, under the circle operation $x\circ
y=x+y+xy$. We proved in \cite{riley-wilson:engel} 
that every finitely generated associative ring $R$ satisfying the Engel 
condition is upper Lie nilpotent.
From this result we were able to infer that whenever $R$ 
satisfies an Engel identity then both the associated circle and multiplicative
semigroups of $R$ must satisfy a so-called Morse identity. 

Define sequences $f_n$ and $g_n$ by $f_1(x,y)=xy, g_1(x,y)=yx$, and 
$$f_{n+1}(x,y,x_3,\dots,x_{n+2})=
f_n(x,y,x_3,\dots,x_{n+1})x_{n+2}g_n(x,y,x_3,\dots,x_{n+1}),$$
$$g_{n+1}(x,y,x_3,\dots,x_{n+2})=
g_n(x,y,x_3,\dots,x_{n+1})x_{n+2}f_n(x,y,x_3,\dots,x_{n+1}),$$
for all
$n\geq 1$. 
The $n$th {\em Malcev identity} (\cite{malcev}) is the semigroup
identity $$f_n(x,y,x_3,\dots,x_{n+1})=g_n(x,y,x_3,\dots,x_{n+1}),$$ 
while the $n$th {\em Morse identity} $u_n(x,y)=v_n(x,y)$ (\cite{morse}) is the $n$th
Malcev identity with $x_3=\cdots=x_{n+1}=1$. 

Consequently, neither $(R,\cdot)$ nor $(R,\circ)$ can contain a free 
subsemigroup if $R$ satisfies an Engel identity. 

The problem of characterizing finitely generated groups satisfying an arbitrary
semigroup identity has been studied by several authors (see, for example,
\cite{lewin-lewin}, \cite{semple-shalev:combinatorial} and \cite{shalev:combinatorial}).
Because this class of groups contains the Burnside groups, this problem 
is highly nontrivial --- 
especially in light of the recent construction by Olshanskii and Storozhev of
a $2$-generated group satisfying a semigroup identity that is not a periodic
extension of a locally soluble group (\cite{olshanskii-storozhev}).

Algebras over fields of characteristic zero which satisfy a circle semigroup law, and a more general semigroup condition called collapsibility, were studied previously by the
first author in \cite{riley:collapsing}.  In sharp contrast to the combinatorial methods employed in  this paper, the techniques used in \cite{riley:collapsing} rest heavily on deep structure theorems from both group and ring theory. In this article we  study associative algebras that satisfy an  arbitrary semigroup identity. In fact, we obtain a partial converse to our result in \cite{riley-wilson:engel}.

Throughout the remainder of this paper, $K$ will denote an infinite commutative 
domain and $R$ an associative $K$-algebra on which the action of $K$ is 
torsion-free (this occurs, for example, when $K$ is an infinite field).  
All identical relations in algebraic
objects will be assumed to be nontrivial unless otherwise stated.
A semigroup $S$ satisfies an identity if and
only if there are distinct words $u,v$ in the free semigroup on
$$X=\{x=x_1,y=x_2,x_3,x_4,\dots\}$$
 so that $u=v$ in $S$. 
The semigroup identity is {\em left reduced} if the first letters of
$u$ and $v$ are different, {\em right reduced} if the last letters of $u$ and $v$ are
different and simply {\em reduced} if it is both left and right reduced. In other
words, $u=v$ is reduced if and only if $uv^{-1}$ and $v^{-1}u$ are reduced words in the
free group on $X$. If $(R,\cdot)$ (respectively $(R,\circ)$)  satisfies an identity we
often say that $R$ satisfies a {\em semigroup identity} (respectively, a {\em circle
semigroup identity}).  Clearly each of these corresponds to a polynomial identity in
$R$. A generalization of a multiplicative semigroup identity  in $R$ is  a 
{\em binomial identity}, a polynomial identity of the form $\alpha_1 u_1+\alpha_2
u_2=0$, where $u_1, u_2$ are monomials and $\alpha_1, \alpha_2\in K$. 
The various types of reduced binomial identities are defined in the obvious way.

V.~Tasi\' c and the first author proved in \cite{riley-tasic:malcev} that
$R$ is Lie nilpotent of class at most $n$ if and only if $(R,\circ)$ satisfies the 
$n$th Malcev identity. 
The main result in the present article further demonstrates the close relationship 
between the Lie structure of $R$ and semigroup properties of $R$:

\begin{theorem}
\label{reduced} Let $R$ be a $K$-algebra. Then the following statements are
equivalent.
\begin{enumerate}
\romenumi

\item $R$ satisfies a circle semigroup identity\textrm{;}

\item $R$ satisfies a reduced semigroup identity\textrm{;}

\item $R$ satisfies a reduced binomial identity\textrm{;}

\item $R$ satisfies an identity of the form $\sum_{i=0}^n \alpha_i y^i x y^{n-i}=0$, 
$\alpha_i\in K, \alpha_0\neq 0, \alpha_n\neq 0$\textrm{;}

\item $R$ satisfies an Engel identity\textrm{;} and,

\item $(R,\circ)$ satisfies a Morse identity.

\end{enumerate}

\end{theorem}

Furthermore, for any two  conditions A, B from (i)--(vi), our proof gives 
(sometimes theoretical) bounds for
the degree of the identity in B in terms of the degree of the identity in A. In
particular, these bounds do not depend on $R$, $K$ or the characteristic of $K$. Notice,
too, that since every finite semigroup (in particular $(R,\circ)$ where $R$ is a finite
ring) satisfies an identity, some hypothesis on the coefficient ring $K$ is required. The
following example demonstrates that the distinction between reduced and arbitrary
multiplicative semigroup identities is also necessary. 

\begin{example}
\label{example}
Let $R$ be the subalgebra of the matrix algebra $M_2(K)$ spanned by the matrix units
$e_{11}$ and $e_{12}$. Then
$[R,R]\subseteq Ke_{12}$, so $R$ satisfies the semigroup identity
$[x,y]z=xyz-yxz=0$. $R$ does not satisfy any Engel
identity, since $[e_{11},e_{12}]=e_{12}$. Thus, by Theorem~\ref{reduced}, $R$ does not
satisfy any reduced semigroup identity, nor any circle semigroup identity. 
\end{example}

\begin{theorem}
\label{general}
Let $R$ be a $K$-algebra where $\ch K=p>0$. Then the following statements are
equivalent.
\begin{enumerate}
\romenumi

\item $R$ satisfies a  semigroup identity\textrm{;}

\item $R$ satisfies a binomial identity\textrm{;}

\item $R$ satisfies an identity of the form $\sum_{i=0}^n \alpha_i y^i x y^{n-i}=0$,
$\alpha_i\in K$\textrm{;} and,

\item $R$ satisfies an identity of the form $y^m e_m y^m=0$. 

\end{enumerate}

\end{theorem}

We remark that the characteristic zero analogue of Theorem~\ref{general} is stated
in \cite{golubchik-mikhalev:semigroup}; however, their result corresponding to our
implication (iv) $\Rightarrow$ (i) is not proved and does not seem
obvious to the present authors.

The fact that $R$ is non-unital is essential to  Example~\ref{example}, as indicated by
the following proposition.
 
\begin{proposition}
\label{unital}
Let $R$ be a unital $K$-algebra. If $R$ satisfies a semigroup identity then
$R$ satisfies the corresponding reduced semigroup identity.
\end{proposition}

\begin{theorem}
\label{upper}
There exists a function $f$, depending only on 
natural numbers $d$ and $n$, such that if a $K$-algebra $R$ satisfies a 
circle semigroup identity of degree $n$ and $R$ is generated over $K$ by
$d$ elements then $R$ is upper Lie nilpotent of index at most $f(d,n)$.
\end{theorem}

%The proofs of the above mentioned results can be found in
%Section~\ref{proofs}.
\section{Semigroup identities}

Our hypotheses on $K$ were chosen to imply, by the usual Vandermonde determinant
argument,  that every homogeneous component of a polynomial identity for $R$ is also a
polynomial identity for $R$ (see \cite[6.4.14]{rowen:rt}). 
We shall use this key fact freely, without explicit mention.

By a {\em partial linear identity} we shall mean an identity of the form 
$$\sum_{i=0}^n \alpha_iy^i x y^{n-i}=0,$$ 
with $\alpha_i\in K$. Such an identity  will be
called {\em left reduced} if $\alpha_0\neq 0$, {\em right reduced} if $\alpha_n \neq
0$ and {\em reduced} if it is both left and right reduced.

\begin{proposition}
\label{linear}Let $R$ be a $K$-algebra.

\begin{enumerate}
\romenumi
\item If a semigroup $S$ satisfies an identity in
$x, y,x_3,\dots$ which is left reduced, right reduced or reduced, then $S$ satisfies an
 identity, of the same type, in $x$ and $y$ only.

\item If $R$ satisfies a binomial identity then $R$ is bounded nil or $R$ satisfies a
semigroup identity. 

\item If $R$ satisfies a  binomial identity which is left reduced, right reduced
or reduced then $R$ satisfies a  partial linear identity of the same type.

\item If $R$ satisfies the identity $y^n=0$  then $R$ satisfies $e_{2n-1}=0$.
\end{enumerate}
\end{proposition}

\begin{proof}
Suppose without loss of generality that our left reduced identity 
has the form
$$x x_{i_1}\cdots x_{i_m}=y x_{j_1}\cdots x_{j_n}.$$
Recall that we identify $x=x_1$ and $y=x_2$.
Substituting $x_i= xy^i, i\geq 3$, we obtain a left reduced identity in $x$ and $y$
only. If the original identity were right reduced as well then $x_{i_m}\neq x_{j_n}$. 
Thus, by an appropriate permutation of the variables, 
we obtain an equivalent identity of the form 
$$x_{k_1}\cdots x_{k_m}x= x_{l_1}\cdots x_{l_n}y.$$ 
Substituting $x_i= xy^i, i\geq 3$ into this identity and
then concatenating on the right with the 2-variable left reduced identity  
yields the 2-variable reduced identity:
$$x x_{i_1}\cdots x_{i_m}x_{k_1}\cdots x_{k_m}x=
y x_{j_1}\cdots x_{j_n}x_{l_1}\cdots x_{l_n}y.$$  
This, and symmetry, proves (i).

Next, given a binomial identity $\alpha_1 u_1+\alpha_2 u_2=0$ holding in $R$, 
set all variables equal, to $y$ say. If the identity is not homogeneous
then separating components shows that $R$ is bounded nil. 
On the other hand, if it is homogeneous then $(\alpha_1 +\alpha_2) y^n=0$ 
for some $n$, so that either $R$ is bounded nil or $\alpha_1=-\alpha_2$, in
which case $u_1-u_2=0$ holds in $R$.  This  proves (ii).

In order to prove (iii), suppose that $R$ satisfies a given binomial identity
and observe from (i) and (ii) that either $R$ is
bounded nil, in which case $R$ satisfies a partial linear identity by (iv) below, 
or $R$ satisfies a semigroup identity of the form $u(x,y)-v(x,y)=0$.
Thus we may assume that $R$ is not bounded nil, and hence that the semigroup identity 
is homogeneous.  We assert that the homogeneous component of degree $1$ in $x$
of the identity $u(x+y,y) - v(x+y,y)=0$ is a partial linear identity.  
To see why it is nontrivial, write $u= au', v= av'$, where
$a$ has length $m$ and $u'=v'$  is a left reduced equation.  
If (as we may assume without loss of generality)
$u'$ starts with $x$ and $v'$ with $y$ then in the expansion of $u(x+y,y)$ there is
precisely one  monomial starting with $y^m x$, whereas no monomial in the
expansion of $v(x+y,y)$ begins with $y^m x$. 
This, and symmetry, yields (iii).

To prove the well-known fact (iv), let $l, r$ denote respectively the $K$-linear
operators of left and right multiplication by $y$. Then, since $l$ and $r$ commute,

\begin{equation}
\label{engel}
e_m=(r-l)^m(x)=\sum_{i=0}^m (-1)^i \binom{m}{i}l^ir^{m-i}(x)=
\sum_{i=0}^m (-1)^i \binom{m}{i}y^ixy^{m-i}.
\end{equation}

Thus if $m=2n-1$ and $R$ satisfies $y^n=0$ then every term in the sum on the right is
zero. 

\end{proof}

Proposition \ref{unital} is a consequence of the following result.

\begin{proposition}
\label{circle}
Let $R$ be a $K$-algebra.
\begin{enumerate}
\romenumi

\item If $(R,\circ)$ satisfies a semigroup identity then $(R,\cdot)$ satisfies the same
identity.
\item If $(R,\circ)$ satisfies a semigroup identity then $(R,\circ)$
satisfies the corresponding reduced identity. 
\item If $R$ is unital then $(R,\circ)\cong (R,\cdot)$.
\end{enumerate}
\end{proposition}

\begin{proof}
Let $S$ be the unital hull of $R$, that is, $S=R$ if $R$ is unital and $S= K1\oplus R$
if $R$ is nonunital. The map $\iota\colon r\mapsto 1+r$ is an injective 
semigroup map from $(R,\circ)$ into $(S,\cdot)$ which is onto if (and only if)
$R=S$. This proves (iii). The image under $\iota$ of an identity in $(R,\circ)$ is an
identity in $(1+R,\cdot)\subseteq (S,\cdot)$. Only the bottom degree homogeneous component
of this identity involves 1 and the other homogeneous components yield identities in
$(R,\cdot)$. The highest degree component is precisely the original identity, yielding
(i).

Assume  that $u(x,y)=v(x,y)$ is an identity for $(R,\circ)$ of degree $n$. Write 
$u=au'b, v=av'b$  where $u'=v'$ is a reduced equation. 
We  show that  $u'=v'$ also holds in $(R,\circ)$.  
It suffices, by symmetry and by induction on the maximum length of $a$ and
$b$,  to prove this in the case when $a=x$ and $b$ is empty. 
The identity $xu'(x,y)=xv'(x,y)$ in $(R,\circ)$ is equivalent to the polynomial identity
$$(1+x)u'(1+x,1+y)-(1+x)v'(1+x,1+y)=0$$ 
in $R$. Let $m$ be an even integer with $m\geq n+1$. Then  multiplying the last identity
on the left by $1-x+x^2-\cdots +x^m$ yields the polynomial identity
$$(1+x^{m+1})u'(1+x,1+y)-(1+x^{m+1})v'(1+x,1+y)=0.$$ 
Separating homogeneous components and using the fact that $x^{m+1}$ has higher
$x$-degree than $u'$ and $v'$, we obtain the polynomial
identity
$$u'(1+x,1+y)-v'(1+x,1+y)=0$$ 
in $R$,  which is equivalent to $u'=v'$ holding in $(R,\circ)$. 
This proves (ii).  
\end{proof}

The following lemma is crucial to our main theorems and is best possible in view of
Example~\ref{example}s. A simpler argument, as in \cite{golubchik-mikhalev:semigroup},
 is available in characteristic zero.  That argument fails in positive characteristic,
where the situation is more delicate.

\begin{lemma}
\label{l1}
Suppose that $R$ satisfies $y^m\alpha y^k=0$ where $\alpha=\sum_i \alpha_i y^i x
y^{n-i}$.
\begin{enumerate}
\romenumi 
\item If $\alpha$ is right reduced  then $R$ satisfies  $y^{m+n} e_n y^k=0$.
\item If $\alpha$ is left reduced then $R$ satisfies  $y^{m} e_n y^{n+k}=0$.
\item If $\alpha$ is reduced then $R$ satisfies $y^me_{3n-1}y^k=0$.
\end{enumerate}
\end{lemma}

\begin{proof}
By symmetry, the proof of (ii) is entirely analogous to that of (i). If the conclusions
of (i) and (ii) hold then the conclusion of (iii) follows from equation~\eqref{engel}:
$$y^me_{3n-1}y^k=y^m\sum_{i=0}^{2n-1}(-1)^i
\binom{2n-1}{i}y^ie_ny^{2n-1-i}y^k=0.$$
Thus it suffices to prove the conclusion of (i).

First assume that $m=k=0$.  Make the substitution $y\mapsto y(y+1)$. Expanding
$\sum_{i=0}^n \alpha_i y^i (y+1)^i x y^{n-i}(y+1)^{n-i}=0$ by the binomial theorem and
separating homogeneous components yields identities $v_0=0,\dots ,v_n=0$ for $R$, where
$v_r$ is homogeneous of degree $n+r$ in $y$. We claim that 

\begin{equation}
\label{amazing}
\sum_{r=0}^n (-1)^r v_r y^{n-r}=\alpha_n y^n e_n.
\end{equation}

To establish equation~\eqref{amazing}, it suffices to  show that the
coefficients of $y^a xy^{2n-a}$ on each side are equal, whenever $0\leq a \leq 2n$.

First note that by equation~\eqref{engel}, the coefficient of $y^a x y^{2n-a}$ in 
$y^n e_n$ is $(-1)^{a-n} \binom{n}{a-n}$ if $a\geq n$ and $0$ otherwise. With the
usual convention on binomial coefficients, the expression $(-1)^{a-n} \binom{n}{a-n}$ is
valid for all $a$. Using the same convention we may sum over all  values of any
index occurring. 

Now we calculate the coefficient of $y^axy^{2n-a}$ in $v_ry^{n-r}$, or, what is the
same, the coefficient of $y^a x y^{n+r-a}$ in $v_r$. The binomial theorem expansion
above shows that the coefficient of $y^s x y^t$  is precisely
$\sum_{i+j=n}\alpha_i\binom{i}{s-i}\binom{j}{t-j}$. Putting $s=a$ and
$t=r+n-a$, we obtain the desired coefficient as 
$\sum_i\alpha_i \binom{i}{a-i}\binom{n-i}{r-(a-i)}$.

It follows that 

\begin{eqnarray*}
&&\sum_r (-1)^r \sum_i \alpha_i \binom{i}{a-i}\binom{n-i}{r-(a-i)}\\
&=&\sum_i \alpha_i \binom{i}{a-i} \sum_r (-1)^r \binom{n-i}{r-(a-i)}\\
&=&\sum_i \alpha_i \binom{i}{a-i}(-1)^{a-i}\sum_s (-1)^s\binom{n-i}{s}\\
&=&(-1)^{a-n}\binom{n}{a-n}\alpha_n
\end{eqnarray*}
 since the inner sum has the value zero unless $n-i=0$, and $1$ otherwise. This proves
(i) in the case $m=k=0$.

In the general case, where $m$ and $k$ are not necessarily zero, the substitution
$y\mapsto y(y+1)$ into the original identity yields an identity

\begin{equation}
\label{unreduced}
\sum_{r+s+t\leq m+n+k} c_{rst}y^s (y^m v_r y^k) y^t=0,
\end{equation}
for some coefficients $c_{rst}\in K$. For $0\leq a\leq n$, consider the homogeneous
component of \eqref{unreduced} of degree $m+n+k+a$ in $y$. The only $v_r$ occurring have
$r\leq a$ and the only term involving $v_a$ is precisely $y^mv_a y^k$. By induction on
$a$,  $y^mv_r y^k=0$ is an identity in $R$ for all $r<a$ and hence so is $y^m v_a
y^k=0$. We may now proceed exactly as in the special case above and the conclusion
follows. 

\end{proof}

\subsection{Unital algebras}
\label{subs:unital}

In case $R$ is unital, more information can be obtained. Note that
$e_n(x,y)=x(\ad y)^n=x(\ad (y+1))^n=e_n(x,y+1)$. Thus by substituting $y\mapsto y+1$
into the result of (i) or (ii)  in Lemma~\ref{l1} and separating out the component of
degree $n$ in $y$ we obtain $e_n=0$ in $R$.

In the rest of this subsection (which is not essential to the main results of the
paper) we give a characterization (for unital $K$-algebras) of the Engel identities.

For each $m\geq 0$, let $W_m$ be the $K$-submodule of $K\langle x, y \rangle$
with basis all monomials $y^i xy^j$ such that $i+j=m$, and let $V_n=\sum_{m=0}^n W_m$
and $V=\sum_{n\geq 0} V_n$. Note that $W_0$ is spanned by the monomial $x$, and that for
$n\geq 1$, $e_n$ is a reduced element of $W_n$.

Define the difference operator $\Delta$ on $V$ by $\Delta
\alpha(x,y)=\alpha(x,y+1)-\alpha(x,y)$. Note that $\Delta\colon V_n\rightarrow V_{n-1}$,
and that the homogeneous component of degree
$n-1$ in $y$ of $\Delta \alpha$ is simply the Hausdorff derivative
$\partial\alpha/\partial y$ with respect to
$y$ (that is, the unique $K$-derivation of $K\langle x,y\rangle$ sending 
$y$ to $1$ and $x$ to $0$). 

\begin{proposition}
Let $R$ be a unital $K$-algebra, and $\alpha\in W_n$.
\begin{enumerate}\romenumi
\item $\Delta \alpha=0$ if and only if $\alpha$ is a scalar
multiple of $e_n$.
\item If $\ch K=0$ then  $\partial\alpha/\partial y=0$ if and only if $\alpha$ is a
scalar multiple of $e_n$. 

\end{enumerate}

\end{proposition}

\begin{proof}Given  $\alpha(x,y)=\sum_{i=0}^n \alpha_i y^i x y^{n-i}$,
expand $\alpha(x,y+1)$ by the binomial theorem. 
The coefficient of $y^s x y^t$ in $\alpha(x,y+1)$ is given by
\begin{equation}
\label{coeff}
\left[y^s x y^t\right]=
\begin{cases}
\sum_i \alpha_i\binom{i}{s}\binom{n-i}{t},& s+t< n\\
0, &s+t=n.
\end{cases}  
\end{equation}
Now  $\Delta \alpha=0$ if and only if the coefficients of all monomials
$y^i x y^j$,  for $i+j\leq n-1$, are zero. This gives a system of 
linear equations in the $n+1$ unknowns $\alpha_0, \dots, \alpha_n$. We claim that the
 coefficient matrix $M$ has rank exactly $n$. Indeed, by
equation~\eqref{coeff}, the submatrix of rows corresponding to the components of $xy^m,
0\leq m\leq n-1$, has the form 

$$
\left[
\begin{matrix}
* &1 &0 &0&\cdots & 0\\
* &* &1 &0&\cdots & 0\\
\vdots &\vdots &\vdots &\ddots & \vdots &0\\
* &* &\dots& *  &1 &0\\
* &* &*& \dots  &* &1
\end{matrix}
\right]
$$ 
which shows that the rank is at least $n$. However the rank is not $n+1$, since, as
observed above, the coefficient vector $\alpha_i=(-1)^i\binom{n}{i}$ of $e_n$ is in the
kernel of $M$. This proves (i).

To prove (ii), it suffices to show that in characteristic zero, the submatrix of $M$
consisting of all rows corresponding to monomials $y^s x y^t$ with $s+t=n-1$ has rank
$n$. By equation~\eqref{coeff}, this submatrix has the form

$$
\left[
\begin{matrix}
n &1 &0&0&\cdots & 0\\
0&n-1 &2&0 &\cdots & 0\\
0&0&n-2&3&\dots&0\\
\vdots &\vdots &\vdots &\ddots & \ddots &0\\
0&0&0 &\dots & 1   &n
\end{matrix}
\right].
$$ 

Since $\ch K=0$ the submatrix consisting of the first $n$ columns is nonsingular and
(ii) follows.
\end{proof}

\section{Proofs of Theorems}
\label{proofs}

We first prove Theorem~\ref{reduced}. The implication (ii) $\Rightarrow$ (iii) is 
obvious and  (i) $\Rightarrow$ (ii) and (iii) $\Rightarrow$ (iv) follow from
Proposition~\ref{circle} and Proposition~\ref{linear} respectively. By Lemma~\ref{l1},
(iv) and (v) are equivalent.  Suppose then that $R$  satisfies $e_n=0$.  
Let $x,y\in R$. By \cite{riley-wilson:engel}, the subalgebra $T$ of $R$ generated by $x$
and $y$ is Lie nilpotent of class $m$ depending on $n$ only. An easy induction on $m$
shows that $T$, and hence $R$, satisfies the Morse identity, in the circle sense, of
degree $m$. Indeed,
$$u_m-v_m=[u_1,v_1, v_2, \dots ,v_{m-1}].$$ 
This proves (v)
$\Rightarrow$ (vi). The last implication  (vi) $\Rightarrow$(i) is obvious.

We now prove Theorem~\ref{general}. The implication (i) $\Rightarrow$ (ii) is
obvious, (ii) $\Rightarrow$ (iii) follows from Proposition~\ref{linear}, and (iii)
$\Rightarrow$ (iv) can be deduced from Lemma~\ref{l1}. 
If $\ch K=p>0$ then (iv) $\Rightarrow$ (i) since by increasing $m$ if necessary  
we may assume that $m=p^t$, so that
$$y^{p^t}xy^{2p^t}-y^{2p^t}xy^{p^t}=y^{p^t}e_{p^t}y^{p^t}=0.$$ 

Finally, Theorem~\ref{upper} follows from the quantitative form of Theorem~\ref{reduced}
and
\cite[Theorem]{riley-wilson:engel}.
\section{Comments}
In an earlier version of this paper, we asked the following questions about an arbitrary
ring $R$.  These questions arose naturally from the work above, 
and the converses had been shown to hold in \cite{riley-wilson:engel}. 

\begin{itemize}
\item If a ring $R$ satisfies a reduced semigroup identity, does $R$
necessarily satisfy an Engel identity?
\item If a ring $R$ satisfies a reduced circle semigroup identity, does $R$
necessarily satisfy an Engel identity?
\end{itemize}

We are indebted to Ol'ga Paison for showing us that the answer to both is no.  
We now present her example.

 Let $p$ be a prime, let $F$ be a field of order $p^2$ and
let $R$ be the subring of $M_2(F)$ consisting of all elements of the form
$ae_{11}+a^pe_{22}+be_{12}$ for $a,b\in F$. Then $R$ does not satisfy any Engel identity.
To see this, choose $a\in F$ with $a^p\neq a$. Let  $x=ae_{11}+a^pe_{22}$ and
$y=e_{12}$. Then for sufficiently large even integers $s$ we have
$[x^{p^s},y]=(a-a^p)e_{12}\neq 0$.  On the other hand, the only idempotents of $R$ are
$0$ and $1$, and so $R$ satisfies a reduced (circle) semigroup identity by the following
result.

\begin{proposition}
Let $R$ be a finite ring. Then
\begin{enumerate}\romenumi
\item $(R,\cdot)$ and $(R,\circ)$ satisfy an identity of the form
 $x^t\equiv x^{2t}$.
\item If all idempotents of $R$ are central, then $R$ satisfies a reduced 
semigroup identity.
\end{enumerate}
\end{proposition}

\begin{proof}The conclusion of part (i) is true for every  finite semigroup  $S$. First,
every element of $S$ is periodic. Furthermore, every periodic element in a semigroup has
some power which is an idempotent. To see this, note that for a fixed $x\in S$,
$x^m=x^{m+a}$ for some $m,a>0$. This implies that for all $n\geq 1$ and all $s\geq m$,
$x^s=x^{s+na}$. Choose $t_x$ such that $t_x\geq m$ and $a$ divides $t_x$. Then
$(x^{t_x})^2=x^{t_x}$.  The desired global identity follows directly from this equation, 
since $S$ satisfies $x^t\equiv x^{2t}$ with $t=\prod_{x\in S} t_x$.

Now by (i), there is some $t$ for which $x^t$ is an idempotent for each $x\in R$. 
Thus if all idempotents of $R$ are central, $R$ satisfies the identity $x^ty\equiv
yx^t$, yielding (ii). 
\end{proof}

In \cite{golubchik-mikhalev:semigroup} it was shown (using arguments
special to characteristic zero) that the  $K$-algebra $R$ satisfies a partial linear
identity if and only if the algebra of $2\times 2$ upper triangular matrices over $K$ is
not in the variety generated by $R$. Perhaps this is true in all characteristics.

\bigbreak

{\bf Acknowledgments.} The first author received support from NSF-EPSCoR in
Alabama and the University of Alabama Research Advisory Committee. The second author 
is supported by a NZST Postdoctoral Fellowship. This work was done while the first 
author visited the Department of Mathematics at the University of Auckland, and their
hospitality is gratefully acknowledged.

\bibliographystyle{amsalpha}
\bibliography{liestruct} 

\end{document}